
\input amstex.tex
\documentstyle{amsppt}
\magnification1200
\hsize=12.5cm
\vsize=18cm
\hoffset=1cm
\voffset=2cm

\def\DJ{\leavevmode\setbox0=\hbox{D}\kern0pt\rlap
{\kern.04em\raise.188\ht0\hbox{-}}D}
\def\dj{\leavevmode
 \setbox0=\hbox{d}\kern0pt\rlap{\kern.215em\raise.46\ht0\hbox{-}}d}

\baselineskip=13pt
\def\hf{{\textstyle{1\over2}}}
\def\a{\alpha}
\def\d{{\,\roman d}}
\def\e{\varepsilon}
\def\f{\varphi}
\def\G{\Gamma}
\def\k{\kappa}
\def\s{\sigma}

\def\={\;=\;}

\def\zt{\zeta(\hf+it)}

\def\D{\Delta}

\def\R{\Re{\roman e}\,}
\def\I{\Im{\roman m}\,}
\def\z{\zeta}

\def\hf{{\textstyle{1\over2}}}

\def\f{\varphi}

\def\le{\leqslant}
\def\ge{\geqslant}
\font\tenmsb=msbm10
\font\sevenmsb=msbm7
\font\fivemsb=msbm5
\newfam\msbfam
\textfont\msbfam=\tenmsb
\scriptfont\msbfam=\sevenmsb
\scriptscriptfont\msbfam=\fivemsb
\def\Bbb#1{{\fam\msbfam #1}}

\def \NN {\Bbb N}
\def \CC {\Bbb C}
\def \RR {\Bbb R}
\def \ZZ {\Bbb Z}

\font\ff=cmr8

\baselineskip=13pt

\font\teneufm=eufm10
\font\seveneufm=eufm7
\font\fiveeufm=eufm5
\newfam\eufmfam
\textfont\eufmfam=\teneufm
\scriptfont\eufmfam=\seveneufm
\scriptscriptfont\eufmfam=\fiveeufm
\def\mathfrak#1{{\fam\eufmfam\relax#1}}

\font\tenmsb=msbm10
\font\sevenmsb=msbm7
\font\fivemsb=msbm5
\newfam\msbfam
     \textfont\msbfam=\tenmsb
      \scriptfont\msbfam=\sevenmsb
      \scriptscriptfont\msbfam=\fivemsb
\def\Bbb#1{{\fam\msbfam #1}}

\def \NN {\Bbb N}
\def \CC {\Bbb C}

\def \RR {\Bbb R}
\def \ZZ {\Bbb Z}

  \def\rightheadline{{\hfil{\ff
On the Rankin-Selberg zeta-function }\hfil\tenrm\folio}}

  \def\leftheadline{{\tenrm\folio\hfil{\ff
   Aleksandar Ivi\'c }\hfil}}
  \def\emptyheadline{\hfil}
  \headline{\ifnum\pageno=1 \emptyheadline\else
  \ifodd\pageno \rightheadline \else \leftheadline\fi\fi}

\font\ff=cmr8
\font\teneufm=eufm10
\font\seveneufm=eufm7
\font\fiveeufm=eufm5
\newfam\eufmfam
\textfont\eufmfam=\teneufm
\scriptfont\eufmfam=\seveneufm
\scriptscriptfont\eufmfam=\fiveeufm
\def\mathfrak#1{{\fam\eufmfam\relax#1}}

\font\tenmsb=msbm10
\font\sevenmsb=msbm7
\font\fivemsb=msbm5
\newfam\msbfam
\textfont\msbfam=\tenmsb
\scriptfont\msbfam=\sevenmsb
\scriptscriptfont\msbfam=\fivemsb
\def\Bbb#1{{\fam\msbfam #1}}

\def \NN {\Bbb N}
\def \CC {\Bbb C}

\def \RR {\Bbb R}
\def \ZZ {\Bbb Z}

\def\D{\Delta}
\def\a{\alpha}
 \def\e{\varepsilon}
 \def\d{\,{\roman d}}
\topmatter
\title
On the Rankin--Selberg zeta-function
\endtitle
\author
Aleksandar Ivi\'c
\endauthor
\address
Katedra Matematike RGF-a, Universitet u Beogradu,  \DJ u\v sina 7,
11000 Beograd, Serbia
\endaddress
\subjclass 11 N 37, 11 M 06, 44 A 15, 26 A 12
\endsubjclass
\email {\tt  ivic\@rgf.bg.ac.rs, aivic\@matf.bg.ac.rs}
\endemail
\dedicatory
\enddedicatory
\abstract We obtain the approximate functional equation
for the Rankin-Selberg zeta-function on the 1/2-line.
\endabstract
\endtopmatter
\document

\head 1. Introduction
\endhead
Let $\varphi(z)$ be a holomorphic cusp form of weight $\kappa$
with respect to the full modular group $SL(2,\ZZ)$, so that
$$
\f\left({az+b\over cz+d}\right) \;= (cz+d)^\k\f(z)\qquad\bigl(a,b,c,d\in\ZZ,\;
 ad-bc=1\bigr)
$$
when $\I z >0$ and $\lim_{\,\I z\to\infty}\f(z)=0$
(see e.g., R.A. Rankin [12] for basic notions).
We denote by $a(n)$ the $n$-th Fourier
coefficient of $\varphi(z)$ and suppose that $\varphi(z)$ is a
normalized eigenfunction for the {\it Hecke operators} $T(n)$, that is,
$a(1)=1$ and $  T(n)\varphi=a(n)\varphi $ for every $n \in
\NN$ (see  Rankin op. cit. for the definition and properties
of the Hecke operators). The classical example is $a(n) = \tau(n)$,
when $\k=12$. This is the well-known {\it Ramanujan tau-function} defined by
$$
\sum_{n=1}^\infty \tau(n)x^n \=
x{\left\{(1-x)(1-x^2)(1-x^3)\cdots\right\}}^{24}\qquad(\,|x| < 1).
$$

Let $c_n\,(\,\ge0)$ be the convolution function defined by
$$
c_{n}=n^{1-\kappa}\sum_{m^2 \mid n}m^{2(\kappa-1)}
\left|a\Bigl({n\over m^2}\Bigr)\right|^2.\leqno(1.1)
$$
Note that $c_n$ is a multiplicative arithmetic function, namely $c_{mn} = c_mc_n$ when
$(m,n)=1$, since $a(n)$ is multiplicative.

The well-known {\it Rankin-Selberg problem} consists
of the estimation of the error term function
$$
\D(x) \;:=\; \sum_{n\leqslant x}c_n - Cx.\leqno(1.2)
$$
The constant $C \,(>0)$ in (1.2) may be written down explicitly
(see e.g., [8]) as
$$
C \;=\;C(\f)\;=\;\frac{2\pi^2(4\pi)^{\kappa-1}}{\Gamma(\kappa)}
      \iint\limits_{\!\mathfrak F} y^{\kappa-2}|\varphi(z)|^{2}\d x\d y,
$$
the integral being taken over a fundamental domain $\mathfrak F$ of
the group $SL(2,\ZZ)$.
The classical upper
bound for $\D(x)$ (strictly speaking $\D(x) = \D(x;\f)$)
of Rankin and Selberg, obtained independently in their important works [11] and [14]
of 1939, is
$$
\D(x) = O(x^{3/5}).\leqno(1.3)
$$
In fact, this result is one of the longest standing unimproved bounds
of Analytic Number Theory, but the present paper is not concerned
with this problem. Our  object of study is
the so-called {\it Rankin--Selberg zeta-function}
$$
Z(s) \;:= \; \sum_{n=1}^\infty c_n n^{-s},            \leqno(1.4)
$$
which is the generating {\it Dirichlet series} for the sequence $\{c_n\}_{n\ge1}$.
One can define the Rankin--Selberg zeta-function in various degrees of generality;
see e.g., Li and Wu [10] where the authors establish universality properties
of such functions.

\smallskip
Note that series in (1.4) converges absolutely for $ \R s > 1$. Namely from
(1.2) and P. Deligne's estimate $|a(n)|\le n^{(\kappa-1)/2}d(n)$ (see [1]),
where $d(n)\,(\ll_\e n^\e)\,$ is the number of positive divisors of $n$, we have
$$
c_n \;\ll_\e\; n^\e,\leqno(1.5)
$$
providing absolute convergence of $Z(s)$ for $\R s > 1$.
\medskip
Here and later $\e$
denotes arbitrarily small constants, not necessarily the same ones
at each occurrence, while $a = O_\e(b)$ (same as $a \ll_\e b$)
 means that the constant implied by the $O$-symbol depends on $\e$.

For  $\R s \le$ the function
$Z(s)$ is defined by analytic continuation.
It has a simple pole at $s=1$ with residue
$C$ (cf. (1.1)), and is otherwise regular. For every $s\in\CC$
it satisfies the functional equation
$$
\G(s+\k-1)\G(s)Z(s) =
(2\pi)^{4s-2}\G(\k-s)\G(1-s)Z(1-s),\leqno(1.6)
$$
where $\G(s)$ is the {\it gamma-function}. One has the decomposition
$$
Z(s)=\zeta(2s)\sum_{n=1}^{\infty}|a(n)|^2n^{1-\kappa-s},
$$
where $\zeta(s)=\sum_{n=1}^\infty n^{-s}\;(\R s>1)$ is the familiar
{\it Riemann zeta-function}.
This formula is the analytic equivalent of the arithmetic relation (1.1).
In our context it is more important that one also has the decomposition
$$
Z(s) := \sum_{n=1}^\infty c_n n^{-s} = \z(s)\sum_{n=1}^\infty b_n
n^{-s} = \z(s)B(s),\leqno(1.7)
$$
say, where $B(s)$ belongs to the {\it Selberg class} of {\it Dirichlet series}
of degree three.  The coefficients $b_n$ in (1.7) are
multiplicative and satisfy
$$
b_n \ll_\e n^\e.\leqno(1.8)
$$
This follows from
$$
b_n = \sum_{d|n}\mu(d)c_{n/d},
$$
which is a consequence of (1.7), the {\it M\"obius inversion formula} and (1.5).
 Actually the coefficients  $b_n$ are bounded by a log-power (see [13])
in mean square, but this stronger property is not needed here.
For the definition
and basic properties of the Selberg class $\Cal S$ of $L$--functions the reader
is referred to A. Selberg's seminal paper [15] and
the comprehensive survey paper of Kaczorowski--Perelli [9].

In view of (1.8) the series for $B(s)$ converges absolutely for $\R s > 1$,
but $B(s)$ has analytic continuation which is holomorphic
for $\R s >0$. This important fact follows from G. Shimura's work [16] (see also A.
Sankaranarayanan [13]), and it implies that (1.7),
namely $Z(s) = \z(s)B(s)$, holds for $\R s >0$ and not only for $\R s>1$.
The function $B(s)$ is of degree three in $\Cal S$, as its functional
equation (see e.g., A.  Sankaranarayanan [13]) is
$$\eqalign{
B(s)\D_1(s) &= B(1-s)\D_1(1-s),\cr \D_1(s) &=
\pi^{-3s/2}\G(\hf(s+\k-1))\G(\hf(s+\k))\G(\hf(s+\k+1)).\cr}
$$
It is very likely that $B(s)$ is primitive in $\Cal S$, namely that it cannot
be factored non-trivially as $F_1(s)F_2(s)$ with $F_1, F_2\in {\Cal S}$,
but this seems hard to prove. Since $B(s)$ is holomorphic
for $\R s>0$, it would follow that one of the factors, say $F_1(s)$, is $L(s+i\a,\chi)$
for some $\a\in\RR$ and $\chi$ a primitive Dirichlet character. This follows from the fact that
elements of degree one in $\Cal S$ are $\z(s+i\a)$ and $L(s+i\a,\chi)$. However, then
$F_2(s)$ would have degree two in $\Cal S$, but the classification of functions in $\Cal S$
of degree two is a difficult open problem.

\medskip
\head 2. The approximate functional equation for $Z(s)$
\endhead
Approximate functional equations are an important tool in the study of
Dirichlet series $F(s) = \sum_{n\ge1}f(n)n^{-s}$. Their purpose
is to approximate $F(s)$ by {\it Dirichlet polynomials} of the type
$\sum_{n\le x}f(n)n^{-s}$ in a certain region where the series defining
$F(s)$ does not converge absolutely.
 In the case of the powers of $\z(s)$  they were studied e.g., in Chapter
 4 of [5] and [6], and in a more general setting by the author [7].


 Before we state our results, which involve approximations of $Z(s)$ by
 Dirichlet polynomials of the form
 $\sum_{n\le x}c_nn^{-s}$, we need some notation. Let (see (1.6))
 $$
X(s) = {Z(s)\over Z(1-s)} = (2\pi)^{4s-2}\,{\G(\k-s)\G(1-s)\over\G(s+\k-1)\G(s)},
\leqno(2.1)
$$
let $\tau = \tau(t)$ be defined by
$$
\log\tau = -{X'(\hf+it)\over X(\hf+it)}\qquad(t\ge3),\leqno(2.2)
$$
and
$$
\Phi(w) = \Phi(w;s,\tau) := \tau^{w-s}X(w) - X(s)\qquad(\hf\le\s = \R s\le1).
\leqno(2.3)
$$
Then we have

\medskip
THEOREM 1. {\it For $\hf\le\s = \R s\le1, t\ge3, s = \s+it$ we have
$$
\eqalign{
Z(s) &= \sum_{n\le x}c_nn^{-s} + X(s)\sum_{n\le y}c_nn^{s-1}
+ C_1{x^{1-s}\over 1-s} + C_2X(s){y^s\over s}\cr&
+ O_\e\Bigl\{t^\e (x^{-\s}+ hx^{1-\s}) +
t^{2+\e-4\s}(y^{\s-1}+hy^\s)\Bigr\}\cr&
-{1\over2\pi ih^3}\int_{{1\over2}-i\infty}^{{1\over2}+i\infty} Z(1-z)\Phi(z;s,\tau)
y^{s-z}(z-s)^{-4}\left(1-{\roman e}^{-h(s-z)}\right)^3\d z,\cr}
\leqno(2.4)
$$
where $xy = \tau, 1\ll x,y\ll \tau, 0 < h \le 1$ is a parameter to be suitably
chosen, and $C_1, C_2$ are absolute constants.}

\medskip
The restriction $\hf\le\s = \R s\le1$ in Theorem 1 can be removed, and one
can consider the whole range $0\le\s\le 1$. For $0\le\s\le\hf$ this is achieved on replacing
$s$ by $1-s$, interchanging $x$ and $y$, and using $Z(1-s)X(s) = Z(s)$,
together with (2.4) and (3.5) of Lemma 2.

\medskip
The most important case of Theorem 1 is when $s=\hf+it$ lies on the so-called
{\it critical line} $\R s = \hf$. Then we obtain from (2.4)
the following

\medskip
THEOREM 2. {\it For $s=\hf+it,\,t\ge 3,\,xy=\tau, 1 \ll x,y \ll \tau$ we have
$$
\eqalign{
Z(s) &= \sum_{n\le x}c_nn^{-s} + X(s)\sum_{n\le y}c_nn^{s-1}
+ C_1{x^{1-s}\over 1-s} + C_2X(s){y^s\over s}\cr&
+  O_\e\Bigl(t^{\e-11/16}(x^{1/2} +
t^2x^{-1/2})^{3/4}\Bigr) + O_\e(t^{1/2+\mu(1/2)+\e}),
\cr}\leqno(2.5)
$$
where, for $\s\in\RR$,}
$$
\mu(\s) \;:=\; \limsup_{t\to\infty}{\log|\z(\s+it)|\over\log t}.
$$

\medskip
The best known result that $\mu(1/2) \le 32/205 = 0.15609\ldots$ is due to M.N.
Huxley [4]. The famous {\it Lindel\"of hypothesis} is that $\mu(1/2)= 0$
(equivalent to $\mu(\s)=0$ for $\s\ge1/2$), and it makes the second error term in (2.5)
equal to $O_\e(t^{1/2+\e})$.

\smallskip
In general, if one introduces smooth weights in the sums in question, then the ensuing
error terms are substantially improved. This was done e.g., in Chapters 4 of [5] and [6]
and in [7]. From the Theorem of [7] (eqs. (19) and (20) with $\s = \hf, K=4,
t\ge 3,\,xy=\tau, 1 \ll x,y \ll \tau$) we obtain
$$
Z(s) = \sum_{n\le x}\rho(n/x)c_nn^{-s} + X(s)\sum_{n\le y}\rho(n/y)c_nn^{s-1}
+  O_\e(t^\e)\;\;(s = \hf+it).\leqno(2.6)
$$
The smooth function $\rho(x)$ (see Chapter 4 of [6] for an explicit construction)
is defined as follows. Let $b>1$ be a fixed constant
 and $\rho(x)\in C^\infty(0,\,\infty)$,
$$
\rho(x) + \rho(1/x) = 1\quad(\forall x\in \RR),\quad \rho(x) = 0 \quad(x\ge b).
$$
\medskip
There is another aspect of this subject worth mentioning.
One can consider the function
$$
{\Cal Z}(t) \;:=\; Z(\hf+it)X^{-1/2}(\hf+it)\qquad(t\in\RR).\leqno(2.7)
$$
The functional equation for $Z(s)$ in the form $Z(s) = X(s)Z(1-s)$ gives easily
$X(s)X(1-s)=1$, hence
$$\eqalign{
\overline{{\Cal Z}(t)} &= Z(\hf-it)X^{-1/2}(\hf-it)
= Z(\hf+it)X(\hf-it)X^{-1/2}(\hf-it)\cr&
= Z(\hf+it)X^{-1/2}(\hf+it)= {\Cal Z}(t).\cr}
$$
Therefore ${\Cal Z}(t)\in\RR$ when $t\in \RR$. The function ${\Cal Z}(t)$
is the analogue of the classical {\it Hardy's function}
$\zt\chi^{-1/2}(\hf+it)$, $\z(s) = \chi(s)\z(1-s)$, which plays a fundamental
role in the study of the zeros of $\z(s)$ on the {\it critical line}
$\R s = 1/2$. Taking $x = (t/2\pi)^2$ in Theorem 2, we obtain then with the
aid of Lemma 2  the following

\medskip
{\bf Corollary}.
$$
 {\Cal Z}(t)= 2\sum_{n\le (t/2\pi)^2}c_nn^{-1/2}\cos\left(t\log
 \bigl({(t/2\pi)^2\over n}\bigr)
 -2t + (\k-1)\pi\right) + O_\e(t^{{1\over2}+\mu({1\over2})+\e}).\leqno(2.8)
 $$

\medskip
One can compare (2.8) to the analogue for $Z^4(t) = |\zt|^4$, since (4.29) of [5]
may be rewritten as
$$
Z^4(t) = 2\sum_{n\le (t/2\pi)^2}d_4(n)n^{-1/2}\cos\left(t\log
\bigl({(t/2\pi)^2\over n}\bigr)-2t-\hf\pi\right) + O_\e(t^{13/48+\e}),\leqno(2.9)
$$
where $d_4(n) = \sum_{abcd=n}1$ is the divisor function generated by $\z^4(s)$. The reason
why the error term in (2.9) is sharper than the one in (2.8) is because we have much more
information on $\z^4(s)$ than on $Z(s)$.
\medskip
The plan of the paper is as follows. In Section 3 we shall formulate and
prove the lemmas necessary for the proofs. In Section 4 we shall prove
Theorem 1, and in Section 5 we shall prove Theorem 2.
\medskip
\head 3. The necessary lemmas
\endhead

\medskip

{\bf Lemma 1}. {\it We have}
$$
\int_0^X|Z(\hf + it)|\d t \;\ll_\e\;X^{5/4+\e}.\leqno(3.1)
$$

\medskip
{\bf Proof of Lemma 1.} From the decomposition (1.7) and the Cauchy-Schwarz
inequality for integrals we obtain
$$
\int_{X/2}^X|Z(\hf + it)|\d t \le \left(\int_{X/2}^X|\z(\hf + it)|^2\d t
\int_{X/2}^X|B(\hf + it)|^2\d t\right)^{1/2}.\leqno(3.2)
$$
Note that we have the elementary bound (see e.g., Chapter 1 of [5])
$$
\int_0^X|\z(\hf + it)|^2\d t \;\ll\; X\log X, \leqno(3.3)
$$
and that $B(s)$ belongs to the Selberg class of degree three. Therefore $B(s)$
is analogous to $\z^3(s)$, and by following the proof of Theorem 4.4 of [5]
(when $k=3$) it is
seen that $B(s)$ satisfies an analogous approximate functional equation,
with $M \ge (3X)^3/Y,\,X^\e\le t\le X$. Taking $Y = X^{3/2}$
and applying the {\it mean
value theorem for Dirichlet polynomials} (Theorem 5.2 of [5]) we obtain,
in view of (1.8), that
$$
\int_{X/2}^X|B(\hf + it)|^2\d t \;\ll_\e\; X^{3/2+\e}.\leqno(3.4)
$$
The bound in (3.1) follows immediately from (3.2)--(3.4)
if we replace $X$ by $X/2^j\;(j = 1,2\ldots)$ and add the resulting expressions.
The best bound for the integral in (3.1) is, up to `$\e$', $X^{1+\e}$.
This follows e.g., by obvious modifications of the arguments used in the proof
of Theorem 9.5 of [5]. It would
improve the bound in (1.3) to $O_\e(x^{1/2+\e})$.

\medskip

{\bf Lemma 2}. {\it For  $0\le\s\le 1$ fixed, $t\ge 3$,  we have}
$$
X(\s+it) = {\Bigl({t\over2\pi}\Bigr)}^{2-4\s}
\exp\left(4it-4it\log\bigl({t\over2\pi}\bigr)+(1-\k)\pi i\right)\cdot
\left(1 + O\left({1\over t}\right)\right), \leqno(3.5)
$$
{\it where the $O$-term admits an asymptotic expansion in
negative powers of $t$}.

\medskip
{\bf Proof of Lemma 2.} Follows from (2.1) and
the full form of Stirling's formula, namely
$$      \log \Gamma(s+b) = (s+b-\hf)\log s - s +\hf\log2\pi
+ \sum_{j=1}^K {(-1)^jB_{j+1}(b)\over j(j+1)s^j}+ O_\delta \biggl( {1\over|s|^{K+1}}\biggr),
$$
which is valid for $b$ a constant, any fixed integer $K \ge 1$,
$|\arg s|\le \pi - \delta$ for $\delta > 0$, where the points $s=0$ and the
neighbourhoods of the poles of~$\Gamma(s+b)$ are excluded, and the $B_{j}(b)$'s
are  {\it Bernoulli polynomials}; for this see e.g., A. Erd\'elyi et al. [2].

\medskip

{\bf Lemma 3}. {\it Let  $\tau = \tau(t)$ be defined by} (2.2). {\it Then
$$
\tau = \Bigl({t\over2\pi}\Bigr)^4\left(1 + O\Bigl({1\over t^2}\Bigr)\right)
\quad(t\ge3),\leqno(3.6)
$$
where the $O$-term admits an asymptotic expansion in negative powers
of $\,t\,$. If $\Phi(w)$ is defined by} (2.3),
{\it then $\Phi(w)(s-w)^{-2}$ is regular for $\R w\le\hf$ and also for $\R w < \s$
if $\hf < \s\le1$.
Moreover, uniformly in $s$ for $\R w = \hf, t\ge3$ we have}
$$
\Phi(w) \;\ll\;t^{2-4\s}\min\Bigl\{1,\min\Bigl(t^{-1}|w-s|^2\Bigr)\Bigr\}.\leqno(3.7)
$$

\medskip
{\bf Proof of Lemma 3.} The functions $\tau$ and $\Phi$ were introduced, in the
case of $\z^2(s)$, by Hardy and Littlewood in their classical proof [3] of the
approximate functional equation for $\z^2(s)$.
To prove (3.6) recall (see (2.1)) that
$$
X(s) = {Z(s)\over Z(1-s)} = (2\pi)^{4s-2}\,{\G(\k-s)\G(1-s)\over\G(s+\k-1)\G(s)}.
$$
Logarithmic  differentiation gives then
$$
\eqalign{
-{X'(\hf+it)\over X(\hf+it)} & = -4\log(2\pi) + {\G'(\k-\hf-it)\over\G(\k-\hf-it)}
+ {\G'(\hf-it)\over\G(\hf-it)}\cr&
+ {\G'(\k-\hf+it)\over\G(\k-\hf+it)}
+ {\G'(\hf+it)\over\G(\hf+it)}.\cr}
$$
If we use (see (A.35) of [5])
$$
{\G'(s)\over\G(s)} = \log s - {1\over2s} + O\Bigl({1\over|s|^2}\Bigr)\qquad
(\,|\arg s|\le \pi-\delta,\,|s|\ge\delta),
$$
where the $O$-term has an asymptotic expansion in term of negative powers of $s$,
we obtain
$$
\log\tau = -{X'(\hf+it)\over X(\hf+it)} = 4\log t - 4\log(2\pi) +
O\Bigl({1\over t^2}\Bigr)\qquad(t\ge3),
$$
which is equivalent to (3.6).

The only non-trivial case concerning the regularity of $\Phi(w)(s-w)^{-2}$ is
when $w = \hf+iv, s = \hf+it$, and this follows from (3.7). For $w = \hf+iv$ we have
$$
|\Phi(w)| \le \tau^{1/2-\s}|X(\hf + iv)| + |X(\s+it)| \ll t^{2-4\s}
$$
in view of (3.6) and (3.5).

To obtain the other bound in (3.7)
suppose that $|w-s|\ll \sqrt{t}$, which is the relevant range of its validity.
Then, for $w=\hf+iv$, we have $v\asymp t$ and
$$
{\d^2\over\d w^2}X(w) \;\asymp\; {1\over t}\qquad(w = \hf+iv,\;v\asymp t).
$$
Write (2.3) as
$$
\Phi(w) \;=\; \tau^{w-s}X(w)\left(1 - {X(s)\over X(w)}\tau^{s-w}\right)\leqno(3.8)
$$
and note that, by Taylor's formula,
$$
\eqalign{
{X(s)\over X(w)}\tau^{s-w}&= \exp\Bigl(\log X(s)-\log X(w) + (s-w)\log\tau\Bigr)\cr&
= \exp\left((s-w){X'(w)\over X(w)} + O(|s-w|^2t^{-1}) + (s-w)\log\tau\right)\cr&
= \exp\left((s-w){X'(\hf+it)\over X(\hf+it)} + O(|s-w|^2t^{-1}) + (s-w)\log\tau\right)
\cr&
= 1 + O\Bigl(|s-w|^2t^{-1}\Bigr),\cr}
$$
in view of (2.3). If we insert this in (3.8)
we obtain  the second estimate in (3.7)
from (3.5), (3.6) and (3.8).

\head
4. Proof of Theorem 1
\endhead

The idea of proof of Theorem 1 goes back to Hardy-Littlewood [3],
who considered the approximate functional equation for $\z^2(s)$.
R. Wiebelitz [17] generalized their method to deal with $\z^k(s)$ when $k\in\NN, k>2$,
and this was refined in Theorem 4.3 of [5]. In what follows we shall make the
modifications which are necessary in the case of $Z(s)$.
 Let the hypotheses of Theorem 1 hold and set
$$
\eqalign{I=I(s,x)&:= {1\over2\pi i} \int_{2-i\infty}^{2+i\infty}Z(s+w)x^ww^{-4}\d w\cr&
= \sum_{n=1}^\infty c_nn^{-s}\left\{{1\over2\pi i}\int_{2-i\infty}^{2+i\infty}
\left({x\over n}\right)^ww^{-4}\d w\right\}\cr&
= {1\over3!}\sum_{n\le x}c_nn^{-s}\log^3(x/n) := S_x,\cr}
$$
say, where we used the absolute convergence of $Z(s)$ for $\s>1$ and
(A.12) of [5] with $k=4$ (reflecting the fact that $Z(s)$ belongs
to the Selberg class of degree $k=4$). The basic idea is
to use a differencing argument to recover $\sum_{n\le x}c_nn^{-s}$ from
the same sum weighted by $\log^3(x/n)$. To achieve this,
first we move the line of integration
in $I$ to $\R w = -1/4$. In doing this we pass over the poles
$w=0$ and $w=1-s$ of the integrand,
with the respective residues
$$
F_x := \sum_{m=0}^3{Z^{(m)}(s)\over m!(3-m)!}(\log x)^{3-m}
$$
and
$$
Q_x : = {Cx^{1-s}\over(1-s)^4}.
$$
Hence by the residue theorem we obtain
$$
J_0:= {1\over2\pi i}\int_{-1/4-i\infty}^{-1/4+i\infty}Z(s+w)x^ww^{-4}\d w
= I-F_x-Q_x = S_x-F_x-Q_x.\leqno(4.1)
$$
In the integral in (4.1) set $z=s+w$, replace $x$ by $\tau/y$, and
use the functional equation for $Z(s)$ and (2.3) in the form
$$
\tau^{u-s}X(u) = X(s) + \Phi(u;s,\tau),
$$
to obtain
$$
\eqalign{
J_0&= {1\over2\pi i}\int_{-1/4-i\infty}^{-1/4+i\infty}Z(1-z)X(s)y^{s-z}(z-s)^{-4}\d z\cr&
+ {1\over2\pi i}\int_{-1/4-i\infty}^{-1/4+i\infty}Z(1-z)
\Phi(z;s,\tau)y^{s-z}(z-s)^{-4}\d z\cr&
= X(s)J_1 + J_2,\cr}
$$
say. This is the point which explains the
definition of the function $\Phi$ in (2.3). We use again
(A.12) of [5] to deduce that
$$
J_1 = {1\over3!}\sum_{n\le y}c_nn^{s-1}\log^3(x/n) := S_y,
$$
similarly to the notation used in evaluating $I$.
The line of integration in $J_2$ is moved to
$\R z = 1/4$. We pass over the pole $z=0$ of the integrand, picking the residue which is
$$
Q_y := - {Cy^s\over s^4}.
$$
Therefore from (4.1) we obtain
$$
F_x - S_x + Q_x = -X(s)(S_y-Q_y)-J_y\leqno(4.2)
$$
with
$$
J_y := {1\over2\pi i}\int_{1/4-i\infty}^{1/4+i\infty}Z(1-z)
\Phi(z;s,\tau)y^{s-z}(z-s)^{-4}\d z.
$$
In (4.2) we replace $x$ and $y$ by $x{\roman e}^{\nu h}$ and $y{\roman e}^{-\nu h}\;
(0\le \nu\le3)$, respectively, so that the condition $x{\roman e}^{\nu h}\cdot
y{\roman e}^{-\nu h}=\tau$ is preserved. We use (see (4.39) and (4.40) of [5])
$$
\sum_{\nu=0}^m (-1)^\nu {m\choose \nu}\nu^p = m!\qquad(p\in\NN)\leqno(4.3)
$$
when $p=m$, and that the sum equals zero when $p<m$, and the estimate
$$
{\roman e}^z  = \sum_{n=0}^M{z^n\over n!} + O(|z|^{M+1})\qquad(M\ge1,\; a\le \R z \le b),
$$
where $a$ and $b$ are fixed. To distinguish better the sums which will arise in this process
we introduce left indices to obtain from (4.2)
$$
\sum_{\nu=0}^3(-1)^\nu{3\choose\nu}\Bigl({}_\nu F_x- {}_\nu S_x+{}_\nu Q_x
+X(s)({}_\nu S_y-{}_\nu Q_y)+{}_\nu J_y\Bigr)=0,
$$
or abbreviating,
$$
{\bar F}_x - {\bar S}_x + {\bar Q}_x + X(s){\bar S}_y -  X(s){\bar Q}_y + {\bar J}_y=0.
\leqno(4.4)
$$
Each term in (4.4) will be evaluated or estimated separately. We have
$$
{\bar F}_x = \sum_{m=0}^3{Z^{(m)}(s)\over 3!(3-m)!}A_m(x)
$$
with
$$
\eqalign{
A_m(x)&:= \sum_{\nu=0}^3(-1)^\nu{3\choose\nu}(\log x+\nu h)^{3-m}\cr&
= \sum_{r=0}^{3-m}{3-m\choose r}h^r\log^{3-m-r}x
\sum_{\nu=0}^3(-1)^\nu{3\choose\nu}\nu^r=3!h^3
\cr}$$
for $m=0$, and otherwise $A_m(x)=0$, where we used (4.3). Therefore
$$
{\bar F}_x = h^3Z(s),
$$
and this is exactly what is needed for the approximate
functional equation that will follow on dividing (4.4) by $h^3$. Consider next
$$
\eqalign{
{\bar S}_x  &= {1\over3!}\sum_{n\le x}c_nn^{-s}\sum_{\nu=0}^3{3\choose\nu}(-1)^\nu
\Bigl(\nu h+\log(x/n)\Bigr)^3\cr&
+ {1\over3!}\sum_{\nu=0}^3{3\choose\nu}(-1)^\nu\sum_{x<n\le x{\roman e}^{\nu h}}
c_nn^{-s}\Bigl(\nu h+\log(x/n)\Bigr)^3\cr&
= \sum\nolimits_1+ \sum\nolimits_2,\cr}
$$
say. Analogously to the evaluation of ${\bar F}_x $ it follows that
$$
\sum\nolimits_1 = h^3\sum_{n\le x}c_nn^{-s}.
$$
We estimate $ \sum\nolimits_2$ trivially, on using (1.5), to obtain
$$
\eqalign{
\Bigl|\sum\nolimits_2\Bigr| &\le {1\over3!}\sum_{\nu=0}^3{3\choose\nu}(2\nu h)^3x^{-\s}
\sum_{x<n\le x{\roman e}^{3h}}c_n\cr&
\ll_\e h^3x^{-\s}t^\e\Bigl(1+ x({\roman e}^{3h}-1)\bigr) \ll_\e t^\e(h^3x^{-\s}
+ h^4x^{1-\s}).\cr}
$$
In a similar way it follows that
$$
\eqalign{&
-X(s){\bar S}_y  = h^3X(s)\sum_{n\le y}c_nn^{s-1}
+ O_\e\Bigl(h^3|X(\s+it)|\sum_{\nu=0}^3\sum_{y{\roman e}^{-3h}<n\le y}c_nn^{\s-1}\Bigr)
\cr&
= h^3X(s)\sum_{n\le y}c_nn^{s-1} +
O_\e\Bigl(h^3t^{2+\e-4\s}(y^{\s-1}+ hy^{\s})\Bigr).
\cr}
$$
Also
$$
{\bar Q}_x = 3!h^3C{x^{1-s}\over1-s} + O(h^4x^{1-\s}),
\quad X(s){\bar Q}_y = C_2X(s)h^3\,{y^s\over s} + O_\e\Bigl(t^{2+\e-4\s}h^4y^{\s}\Bigr).
$$
Therefore we are left with the evaluation of
$$
{\bar J}_y =
{1\over2\pi i}\int_{1/4-i\infty}^{1/4+i\infty}Z(1-z)\Phi(z;s,\tau)y^{s-z}(z-s)^{-4}
\sum_{\nu=0}^3(-1)^\nu{3\choose\nu}{\roman e}^{-\nu h(s-z)}\d z.
$$
Observing that (3.7) holds and that the function
$$
\sum_{\nu=0}^3(-1)^\nu{3\choose\nu}{\roman e}^{-\nu h(s-z)}
= \Bigl(1- {\roman e}^{- h(s-z)}\Bigr)^3
$$
has a zero of order three at $z=s$, we can move the line of integration in
${\bar J}_y $ to $\R z =\hf$. Hence
$$
{\bar J}_y =
{1\over2\pi i}\int_{1/2-i\infty}^{1/2+i\infty}Z(1-z)\Phi(z;s,\tau)y^{s-z}(z-s)^{-4}
\Bigl(1- {\roman e}^{- h(s-z)}\Bigr)^3\d z.
$$
Therefore we obtain
the assertion of Theorem 1 from (4.4) if we divide the whole expression by $h^3$
and collect the above estimates for the error terms.

\head 5. Proof of Theorem 2
\endhead
\medskip
We set $s=\hf+it, z = \hf+iv$ in (2.4) and write the integral on the right-hand side  as
$$
i\int_{-\infty}^\infty\cdots\d v = i\Biggl(\int_{-\infty}^{t/2} + \int_{t/2}^{2t}
+ \int_{2t}^\infty\Biggr)\cdots\d v
= i\Bigl(I_1+I_2+I_3\Bigr), \leqno(5.1)
$$
say. The integrals $I_1$ and $I_3$ are estimated analogously. The latter is, by trivial estimation
and the first bound in (3.7),
$$
\eqalign{&
\int_{2t}^\infty Z(\hf-iv)\Phi(\hf+iv;s,\tau)
y^{i(t-v)}(t-v)^{-4}\left(1-{\roman e}^{-hi(t-v)}\right)^3\d v\cr&
\ll \int_{2t}^\infty |Z(\hf+iv)|v^{-4}\d v \ll_\e t^{\e-11/4},\cr}\leqno(5.2)
$$
where we used (3.1) of Lemma 1. From (2.4), (5.1) and (5.2) it follows that
$$
\eqalign{
Z(s) &= \sum_{n\le x}c_nn^{-s} + X(s)\sum_{n\le y}c_nn^{s-1}
+ C_1{x^{1-s}\over 1-s} + C_2X(s){y^s\over s}\cr&
+ O_\e\Bigl(1 + t^{\e-11/16}(x^{1/2} + t^2x^{-1/2})^{3/4}\Bigr)
- {1\over2\pi ih^3}I_2,
\cr}\leqno(5.3)
$$
with the choice
$$
h \;=\; t^{-11/16}(x^{1/2}+ t^2x^{-1/2})^{-1/4},
$$
so that $0 < h \le 1$ holds. To estimate $I_2$ we use
$$
\left(1-{\roman e}^{-hi(t-v)}\right)^3 \ll h^3|t-v|^3
$$
and the second bound in (3.7) ($\s=\hf$). This gives, on using
the Cauchy-Schwarz inequality for integrals,
$$
\eqalign{
h^{-3}I_2 &\ll \int_{t/2}^{2t}|Z(\hf+iv)|\min\left({1\over|t-v|},{|t-v|\over v}\right)\d v\cr&
\ll {\left(\int_{t/2}^{2t}|Z(\hf+iv)|^2\d v\right)}^{1/2}
{\Bigl(j_1+j_2+j_3\Bigr)}^{1/2},\cr}\leqno(5.4)
$$
say. By (1.7), (3.4) and the definition of the $\mu$-function we have
$$
\int_{t/2}^{2t}|Z(\hf+iv)|^2\d v = \int_{t/2}^{2t}|\zeta(\hf+iv)|^2|B(\hf+iv)|^2\d v
\ll_\e t^{2\mu(1/2)+3/2+\e}.\leqno(5.5)
$$
We have
$$
j_1 := \int_{t/2}^{t-\sqrt{t}}{\d v\over (t-v)^2} \ll {1\over\sqrt{t}},
$$
and the same bound holds for
$$
j_3 :=\int_{t+\sqrt{t}}^{2t}{\d v\over (t-v)^2}.
$$
We also have
$$
j_2 \;:=\; \int_{t-\sqrt{t}}^{t+\sqrt{t}}{\bigl(t-v\bigr)}^2{\d v\over v^2}
\;\ll\; {1\over\sqrt{t}},
$$
so that from (5.4), (5.5) and the bounds for $j_1, j_2, j_3$ we infer that
$$
h^{-3}I_2 \;\ll_\e\; t^{1/2+\mu(1/2)+\e}.\leqno(5.6)
$$
The assertion of Theorem 2 follows then from (5.3) and (5.6),
since the first error term in (5.3) is absorbed by the right-hand side of
(5.6) because $x^{1/2} \ll t^2$.

\vfill\eject

\Refs
\bigskip

\item{[1]} P. Deligne, La conjecture de Weil,
    Inst. Hautes \'{E}tudes Sci. Publ. Math. {\bf 43} (1974), 273-307.

\item{[2]} A. Erd\'elyi, W. Magnus, F. Oberhettinger and F.G. Tricomi, Higher
Transcendental Functions. Volume I, McGraw-Hill, 1953.

\item{[3]} G.H. Hardy and J.E. Littlewood, The approximate functional equation
for $\z(s)$ and $\z^2(s)$, Proc. London Math. Soc. (2){\bf29}(1929), 81-97.

\item{[4]} M.N. Huxley,
Exponential sums and the Riemann zeta function V,
Proc. London Math. Soc. (3) {\bf 90}(2005), 1-41.

\item{[5]} A. Ivi\'c,  The Riemann Zeta-Function
   John Wiley \& Sons, New York, 1985 (2nd ed. Dover, Mineola, New York, 2003).

   \item{[6]} A. Ivi\'c,  The mean values of the Riemann zeta-function,
LNs {\bf 82}, Tata Inst. of Fundamental Research, Bombay 1991 (distr. by
Springer Verlag, Berlin etc.).

  \item{[7]} A. Ivi\'c, An approximate functional equation for a class of
  Dirichlet series, J.  Analysis (Madras, India) {\bf3} (1995), 241-252.

\item{[8]} A. Ivi\'c, K. Matsumoto and Y. Tanigawa, On Riesz mean
of the coefficients of the Rankin--Selberg series, Math. Proc.
Camb. Phil. Soc. {\bf127}(1999), 117-131.

\item{[9]} A. Kaczorowski and A. Perelli, The Selberg class: a
survey, in ``Number Theory in Progress, Proc. Conf. in honour
of A. Schinzel (K. Gy\"ory et al. eds)", de Gruyter,
Berlin, 1999, pp. 953-992.

\item{[10]} Li, Hongze and Wu, Jie,
The universality of symmetric power $L$-functions and their
Rankin-Selberg $L$-functions, J. Math. Soc. Japan {\bf59}(2007), 371-392.

\item{[11]}R.~A.~Rankin, Contributions to the theory of Ramanujan's
   function $\tau(n)$ and similar arithmetical functions II. The order
   of the Fourier coefficients of integral modular forms,
    Proc. Cambridge Phil. Soc. {\bf 35}(1939), 357-372.

\item{[12]}R.~A.~Rankin, Modular forms and functions, Cambridge Univ. Press,
Cambridge, 1977.

\item{[13]} A. Sankaranarayanan, Fundamental properties of symmetric
square $L$-functions I, Illinois J. Math. {\bf46}(2002), 23-43.

    \item{[14]}A.~Selberg,  Bemerkungen \"{u}ber eine Dirichletsche Reihe,
    die mit der Theorie der Modulformen nahe verbunden ist,
    Arch. Math. Naturvid. {\bf 43}(1940), 47-50.

    \item{[15]} A.~Selberg, Old and new conjectures and results about a class
    of Dirichlet series, in ``Proc. Amalfi Conf. Analytic Number Theory 1989
    (E. Bombieri et al. eds.)",
    University of Salerno, Salerno, 1992, pp. 367--385.

    \item{[16]} G. Shimura, On the holomorphy of certain Dirichlet series,
    Proc. London Math. Soc. {\bf31}(1975), 79-98.

\item{[17]} R. Wiebelitz, \"Uber approximative Funktionalgleichungen der
Potenzen der Riemannschen Zeta-funktion, Math. Nachr. {\bf6}(1951-1952),
263-270.

\vskip1cm
\endRefs

\enddocument

\end